\documentclass[12pt, leqno]{article}
\usepackage{amssymb,amsmath}
\usepackage[pdftex]{color}
\usepackage[pdftex]{graphicx}
\usepackage{latexsym}
\usepackage{graphicx}
\usepackage{epstopdf}
\newtheorem{theorem}{Theorem}

\newtheorem{definition}{Definition}

\newtheorem{remark}{Remark}
\newtheorem{lemma}{Lemma}

\newcommand{\tx}[1]{\mbox{\;{#1}\;}} 

\newcommand{\N}{\mathbb{N}}
\newcommand{\R}{\mathbb{R}^n}
\newcommand{\rz}{\mathbb{R}}
\newcommand{\K}{{\cal{K}}}
\newcommand{\E}{{\cal{E}}}
\newcommand{\F}{{\cal{F}}}
\newcommand{\Div}{\hbox{div}\:}
\newcommand{\p}{\partial}
\numberwithin{equation}{section}
\begin{document}
\title{Shape optimization for an elliptic operator with infinitely many positive and negative eigenvalues}
\pagestyle{myheadings}
\maketitle
\centerline{\scshape Catherine Bandle}
\medskip
{\footnotesize
\centerline{Departement Mathematik und Informatik, Universit\"at Basel,}
\centerline{Spiegelgasse 1,  CH-4051 Basel, Switzerland}
} 
\medskip
\centerline{\scshape Alfred Wagner}
\medskip
{\footnotesize
\centerline{Institut f\"ur Mathematik, RWTH Aachen  }
\centerline{Templergraben 55, D-52062 Aachen, Germany}}
\bigskip

\abstract{}
\bigskip

\bigskip

\centerline {MSC2010:  49K20, 49R05, 15A42, 35J20, 35N25.}
\abstract{The paper deals with an eigenvalue problems possessing infinitely many positive and negative
eigenvalues. Inequalities for the smallest positive and the largest negative eigenvalues, which have the same
properties as the fundamental frequency, are derived. The main question is whether or not the classical isoperimetric
inequalities for the fundamental frequency of membranes hold in this case. The arguments are based on the harmonic transplantation
for the global results and the shape derivatives (domain variations) for nearly circular domain. }
\section{Introduction}
In this paper we study the spectrum of the problem
\begin{eqnarray}
\label{resolvent}\Delta u+ \lambda u= 0\qquad\hbox{in}\:\Omega\qquad\qquad\partial_{\nu} u= \lambda\:\sigma\: u\qquad\hbox{in}\:\partial\Omega,
\end{eqnarray}
where $\sigma, \lambda\in\rz$, and $\Omega\subset\R$ is a bounded domain with smooth boundary. The corresponding Rayleigh quotient is
\begin{eqnarray}\label{ray}
R_{\sigma}(v)=\frac{\int\limits_{\Omega}\vert\nabla v\vert^2\:dx}{\int\limits_{\Omega}v^2\:dx+\sigma\int\limits_{\partial\Omega}v^2\:dS}\qquad v\in H^{1,2}(\Omega)\:.
\end{eqnarray}
For positive $\sigma$ the Rayleigh quotient is positive and the classical theory for  symmetric operators applies. Fran\c{c}ois \cite{Fr06} has shown hat in this case the spectrum consists of countably many eigenvalues, which are bounded from below and tend to infinity. 

Our interest in this paper is directed to the problem with  $\sigma$ negative. It has been studied in \cite{BaBeRe06} and \cite{BR04}. A more general approach is found in
\cite{ErSc65}. It is known that in addition to $\lambda=0$ there exist two sequences  of eigenvalues, one tending to $+\infty$ and the other tending to $-\infty$. The eigenfunctions are complete in $H^1(\Omega)$ except in the {\sl resonance case} $|\Omega| +\sigma |\p \Omega|=0$ where some supplements are required (see \cite{BaBeRe06}).

The smallest positive eigenvalue $\lambda_1(\Omega)$ and the largest negative eigenvalue $\lambda_{-1}(\Omega)$ play the role of fundamental frequencies. Based on the isoperimetric inequalities for the fundamental frequency of the membrane we study the dependence of $\lambda_{\pm 1}(\Omega)$ on some geometrical properties such as the volume and the harmonic radius. 

We first establish inequalities by means of the harmonic transplantation which is appropriate for this type of problems. An interesting question is whether the Rayleigh- Faber -Kahn inequality extends to these eigenvalues. Here only  answers for nearly spherical domains can be provided. The arguments are based on the first and second order shape derivatives. For general domains the answer is still incomplete.
 \medskip

One motivation for studying this problem are dynamical boundary conditions for parabolic equations. A simple version is given by the heat equation
\begin{eqnarray*}
\partial_{t}u-\Delta u &=& 0\qquad\hbox{in}\:(0,\infty)\times\Omega\\
\sigma\:\partial_{t}u+\partial_{\nu} u&=& 0\qquad\hbox{in}\:(0,\infty)\times\partial\Omega\\
u(x,0)&=&u_{0}(x)\qquad\hbox{in}\:\Omega.
\end{eqnarray*}
Such problem are well studied (see e.g. \cite{E89} - \cite{E93}, \cite{Hi89},\cite{VV08},\cite{VV09}). It is known that they are well posed for positive $\sigma$ in the space $C(0,T, H^{1,2}(\Omega))$, in the sense of Hadamard, and that there exists a smooth solution globally in time, whereas this is not the case if $\sigma<0$. That is, there is no continuous dependence on the initial conditions (except in dimension one). 
\medskip

The paper is organized as follows. First we present the eigenvalue problem and quote some known results. Then we derive inequalities by means of the method of harmonic transplantation. In the last part we compute the first and the second domain variations of the fundamental eigenvalues and derive some inequalities and monotonicity properties of nearly circular domains.

\section{The eigenvalue problem, known results}
Let $\Omega\subset\R$ be a bounded domain with smooth boundary. Assume $\sigma>0$. In that case \eqref{ray} gives a non-negative functional $R_{\sigma}$ on the Hilbert space $H^{1,2}(\Omega)$. Its minimum is equal to zero and is achieved by any constant function. If we minimize $R_{\sigma}$ over the set
\begin{eqnarray*}
\left\{u\in H^{1,2}(\Omega) : \int\limits_{\Omega}u\:dx=0\right\}
\end{eqnarray*}
the direct method of the calculus of variations gives us a unique minimizer which is non-negative and solves \eqref{resolvent}. By Harnack's inequality it is positive in $\Omega$. G. Fran\c{c}ois in \cite{Fr06} showed that there exist countably many eigenvalues 
$(\lambda_{n})_n$ which tend to $\infty$, for which \eqref{resolvent} admits solutions. 
\medskip

The case $\sigma <0$ was considered in \cite{Ba08} and \cite{BaBeRe06} and will be under 
consideration from now on. For $u,v\in H^{1,2}(\Omega)$ let
\begin{eqnarray*}
a(u,v):=\int\limits_{\Omega}u\:v\:dx+\sigma\int\limits_{\partial\Omega}u\:v\:dS
\end{eqnarray*}
be an inner product on $L^2(\Omega)\oplus L^2(\partial\Omega)$. We define
\begin{eqnarray*}
\K:=\left\{u\in H^{1,2}(\Omega)\::\:\int\limits_{\Omega}\vert\nabla u\vert^2\:dx= 1,\:a(u,1)=0\right\}.
\end{eqnarray*}
From \eqref{ray} we then have 
\begin{eqnarray*}
R_{\sigma}(u)=\frac{1}{a(u,u)}\qquad\hbox{for all}\: u\in\K.
\end{eqnarray*}
\medskip
 
In \cite{BR04} and \cite{BaBeRe06} the authors showed the existence of two infinite sequences of eigenvalues. One sequence consists of negative eigenvalues  $(\lambda_{-n})_n$ and the other of positive eigenvalues $(\lambda_{n})_n$. The corresponding eigenfunctions $(u_{\pm n})_n\in\K$ solve \eqref{resolvent}. Moreover
\begin{eqnarray*}
\lim_{n\to\infty}\lambda_{-n}=-\infty\qquad\hbox{and}\qquad\lim_{n\to\infty}\lambda_{n}=\infty.
\end{eqnarray*}
The eigenvalues are ordered as
\begin{eqnarray*}
\hdots\leq\lambda_{-n}\leq\hdots\leq\lambda_{-1}<\lambda_0=0<\lambda_{1}\leq\lambda_{2}\leq\hdots\leq\lambda_n\leq\ldots.
\end{eqnarray*}
Note that $\lambda_{\pm n}=\lambda_{\pm n}(\sigma)$ for $n\in\N$. 
\medskip

In \cite{BaBeRe06} it was also shown, that the quantity
\begin{eqnarray*}
\sigma=\sigma_{0}:=-\frac{\vert\Omega\vert}{\vert\partial\Omega\vert}.
\end{eqnarray*}
plays an important role. If $\sigma\ne \sigma_0$ the following characterization of $\lambda_{\pm1}$ holds :
\begin{eqnarray}\label{eigenchar}
\lambda_{1}(\sigma)=\frac{1}{\sup_{u\in\K} a(u,u)}>0,\qquad\qquad
\lambda_{-1}(\sigma)=\frac{1}{\inf_{u\in\K} a(u,u)}<0.
\end{eqnarray}
Moreover
\begin{eqnarray*}
\lim_{\sigma\to\sigma_{0}}\lambda_{\pm 1}=0.
\end{eqnarray*}

For $n\ne 1$ (and  $\sigma\ne \sigma_0$) the eigenvalues $\lambda_{\pm n}$ have a variational 
characterization as well. Indeed, let $\lambda_{\pm i}$, $i=1,\hdots,k$ be the first $k$ eigenvalues, 
counted with their multiplicities. Let $u_i$ be the corresponding 
eigenfunctions. Then we define
\begin{eqnarray*}
\K_{k}:=\left\{u\in H^{1,2}(\Omega)\::\:\int\limits_{\Omega}\vert\nabla u\vert^2\:dx= 1,\:a(u,1)=0, a(u,u_{i})=0,\:i=1\hdots k\right\}.
\end{eqnarray*}
We get the characterization
\begin{eqnarray*}
\lambda_{k+1}(\sigma)=\frac{1}{\sup_{u\in\K_{k}} a(u,u)}>0,\qquad\qquad
\lambda_{-k-1}(\sigma)=\frac{1}{\inf_{u\in\K_k} a(u,u)}<0.
\end{eqnarray*}
\begin{remark}\label{sigmacon}
It is interesting to note, that in the case $\sigma\ne\sigma_0$ the constraint 
\begin{eqnarray*}
a(u_{\pm 1},1)=0
\end{eqnarray*} 
is satisfied automatically. This can be seen by the following considerations. Let $u_{\pm 1}$ be an 
eigenfunction corresponding to $\lambda_{\pm 1}$. Then necessarily
\begin{eqnarray*}
\int\limits_{\Omega}u_{\pm 1}v\:dx+\sigma\int\limits_{\partial\Omega}u_{\pm 1}v\:dS-
\mu_{\pm 1}\int\limits_{\Omega}\nabla u_{\pm 1}\cdot\nabla v\:dx-\tilde{\mu}_{\pm 1}a(v,1)=0
\end{eqnarray*}
for all $v\in H^{1,2}(\Omega)$. The real numbers $\mu_{\pm 1}$ and $\tilde{\mu}_{\pm 1}$ denote the Lagrange parameters for the constraints. The special choice $v=u_{\pm 1}\in\K$ gives
\begin{eqnarray*}
\mu_{\pm 1}=\frac{1}{\lambda_{\pm 1}}.
\end{eqnarray*}
For $v=1$ we get
\begin{eqnarray*}
\tilde{\mu}_{\pm 1}\left(\vert\Omega\vert+\sigma\vert\partial\Omega\vert\right)=0.
\end{eqnarray*}
For $\sigma\ne \sigma_{0}(\Omega)$ this implies $\tilde{\mu}_{\pm 1}=0$. Thus in this case the constraint $a(u_{\pm 1},1)=0$ is automatically satisfied. 
\end{remark}
\begin{remark}\label{sigmacon2}
Let $\lambda_1$ be given by \eqref{eigenchar}, i.e.
\begin{eqnarray*}
\frac{1}{\lambda_{1}(\sigma)}=\sup_{u\in\K} a(u,u).
\end{eqnarray*}
Assume $\sigma<\sigma_{0}$ and let $u\in H^{1,2}(\Omega)$ satisfy $\int\limits_{\Omega}\vert\nabla u\vert^2\:dx=1$ only. Then $u$ it is not admissible, since $a(u,1)$ may be different from zero. Let $c\in\rz$ be chosen as
\begin{eqnarray*}
c:=-\:\frac{a(u,1)}{a(1,1)}.
\end{eqnarray*}
With this choice we have $u+c\in \K$ and
\begin{eqnarray*}
a(u+c,u+c)=a(u,u)+2c\: a(u,1)+c^2 a(1,1)=a(u,u)-\:\frac{a(u,1)^2}{a(1,1)}\geq a(u,u).
\end{eqnarray*}
By our assumptions on $\sigma$ we get
\begin{eqnarray}\label{bawa1}
\frac{1}{\lambda_{1}(\sigma)}\geq a(u+c,u+c)\geq a(u,u).
\end{eqnarray}
In the same way we prove 
\begin{eqnarray}\label{bawa2}
\frac{1}{\lambda_{-1}(\sigma)}\leq a(u,u)
\end{eqnarray}
for $0>\sigma>\sigma_0$.
\end{remark}
In \cite{BaBeRe06} (Theorem 14, Theorem 21 and Corollary 22) the following result was proved.
\begin{lemma}\label{sigma0}
\begin{itemize}
\item[(i)] If $\sigma<\sigma_0(\Omega)<0$, then $\lambda_{1}(\sigma)$ is simple and the corresponding eigenfunction $u_{1}$ is of constant sign.
\item[(ii)] If $\sigma_0(\Omega)<\sigma<0$, then $\lambda_{-1}(\sigma)$ is simple and the corresponding eigenfunction $u_{-1}$ is of constant sign.
\item[(iii)] If $\sigma=\sigma_0(\Omega)$ then both $u_{1}$ and $u_{-1}$ change sign. In particular 
\begin{eqnarray*}\lambda_{1}(\sigma_0)=\lambda_{-1}(\sigma_0)=0.
\end{eqnarray*}
\end{itemize}
\end{lemma}
\begin{lemma}\label{sigmamon}
The eigenvalues $\lambda_{\pm 1}(\sigma)$ are monotonically decreasing functions of $\sigma$.
\end{lemma}
{\bf{Proof}} We assume $\sigma_1>\sigma_2$. We distinguish two cases.
\newline
Case 1. From the characterization of $\lambda_{1}$ we get
\begin{eqnarray*}
\frac{1}{\lambda_{1}(\sigma_1)}
\geq
\int\limits_{\Omega}u^2\:dx+\sigma_1\int\limits_{\partial\Omega}u^2\:dS
\geq
\int\limits_{\Omega}u^2\:dx+\sigma_2\int\limits_{\partial\Omega}u^2\:dS.
\end{eqnarray*}
For $u$ we choose the eigenfunction of $\lambda_{1}(\sigma_2)$ and we obtain
\begin{eqnarray*}
\frac{1}{\lambda_{1}(\sigma_1)}
\geq
\frac{1}{\lambda_{1}(\sigma_2)}.
\end{eqnarray*}
This gives $\lambda_{1}(\sigma_1)\leq \lambda_{1}(\sigma_2)$.
\newline
Case 2. From the characterization of $\lambda_{-1}$ we get
\begin{eqnarray*}
\frac{1}{\lambda_{-1}(\sigma_2)}
\leq
\int\limits_{\Omega}u^2\:dx+\sigma_2\int\limits_{\partial\Omega}u^2\:dS
\leq
\int\limits_{\Omega}u^2\:dx+\sigma_1\int\limits_{\partial\Omega}u^2\:dS.
\end{eqnarray*}
In this case we choose $u$ as the eigenfunction of $\lambda_{-1}(\sigma_1)$ and we obtain
\begin{eqnarray*}
\frac{1}{\lambda_{-1}(\sigma_2)}
\leq
\frac{1}{\lambda_{-1}(\sigma_1)}.
\end{eqnarray*}
Since  $\lambda_{-1}(\sigma)<0$ we have $\lambda_{-1}(\sigma_1)\leq \lambda_{-1}(\sigma_2)$.
\hfill $\square$
\begin{remark}\label{asymp}
In \cite{BaBeRe06} the authors also studied the smoothness and asymptotic behaviour of the map
\begin{eqnarray*}
\sigma\to\lambda_{\pm 1}(\sigma).
\end{eqnarray*}
They proved that
\begin{align}
\lambda(\sigma)
= 
\begin{cases} 
\lambda_{1}(\sigma)&  \tx{if} \sigma<\sigma_{0}\\
0&  \tx{if} \sigma=\sigma_{0}\\
\lambda_{-1}(\sigma)&  \tx{if} 0>\sigma>\sigma_{0}
\end{cases}
\end{align}
is a smooth curve with the following asymptotics:
\begin{eqnarray*}
\lim_{\sigma\to-\infty}\lambda(\sigma)=\mu_{D}\qquad\hbox{and}\qquad\lim_{\sigma\to 0}\lambda(\sigma)=-\infty,
\end{eqnarray*}
where $\mu_D$ is the first Dirichlet eigenvalue for the Laplacian.
\end{remark}
We are interested in the domain dependence of $\lambda_{\pm1}$. Thus we will write
$\lambda_{1}=\lambda_{1}(\Omega)$ and $\lambda_{-1}=\lambda_{-1}(\Omega)$. 
Note that $\sigma_0=\sigma_0(\Omega)$ depends on $\Omega$ as well. Moreover for domains of given volume and for a ball $B_R$ with the same volume the isoperimetric inequality gives
\begin{eqnarray}\label{sigmaiso}
\sigma_0(\Omega)
=
-\frac{\vert\Omega\vert}{\vert\partial\Omega\vert}
=
-\frac{\vert B_R\vert}{\vert\partial\Omega\vert}
\geq
-\frac{\vert B_R\vert}{\vert\partial B_R\vert}
=
\sigma_0(B_R).
\end{eqnarray}
In \cite{Ba08} and \cite{BaBeRe06}  the following properties were proved.
\begin{lemma}\label{domaindep}
For some given ball $B_R$ and $\sigma_0(B_R)=-\frac{R}{n}$ let $\lambda_{1}$ and $\lambda_{-1}$ be given by \eqref{eigenchar}. Then the following cases occur.
\begin{itemize}
\item[(i)] Let $B_R$ be a ball such that $\vert B_R\vert=\vert \Omega\vert$. If $\sigma<\sigma_0(B_R)$, then 
\begin{eqnarray}\label{glopt1}
\lambda_{1}(\Omega)\geq \lambda_{1}(B_R).
\end{eqnarray}
\item[(ii)] For any domain $\Omega$ with the same volume as $\vert B_R\vert$, there exist a number $\hat{\sigma}\in(\sigma_0(\Omega),0)$ such that 
\begin{eqnarray}\label{glopt2}
\lambda_{-1}(\Omega)\geq\lambda_{-1}(B_R)
\end{eqnarray} 
whenever $\sigma\in(\sigma_0(\Omega),\hat{\sigma})$.
\end{itemize}
\end{lemma}
\begin{remark}\label{isop}
For (i) we note that the condition $\sigma<\sigma_0(B_R)$ is more restrictive than the condition $\sigma<\sigma_0(\Omega)$ if $\vert\Omega\vert=\vert B_R\vert$. This is a consequence of \eqref{sigmaiso}.
\end{remark}
\section{Harmonic transplantation}
From Section 2 we know that the eigenvalues $\lambda_{1}$ (resp. $\lambda_{-1}$) have a variational characterization for $\sigma\ne\sigma_0(\Omega)$. Moreover for $\sigma<\sigma_0(\Omega)$
\begin{eqnarray*}
0\leq\frac{1}{\lambda_{1}(\Omega)}=\sup\left\{\int\limits_{\Omega}v^2\:dx-\vert\sigma\vert\int\limits_{\partial\Omega}v^2\:dS\::\: v\in \K\right\}
\end{eqnarray*}
has a unique minimizer (of constant sign). The same holds in the case $\sigma_0(\Omega)<\sigma<0$ for
\begin{eqnarray*}
0\geq\frac{1}{\lambda_{-1}(\Omega)}=\inf\left\{\int\limits_{\Omega}v^2\:dx-\vert\sigma\vert\int\limits_{\partial\Omega}v^2\:dS\::\: v\in \K\right\}.
\end{eqnarray*}
We are interested in optimality results for these eigenvalues. They will be obtained by means of the method of harmonic transplantation which was introduced by  Hersch\cite{He69}, (cf. also \cite{BaFl96}). It generalizes the conformal transplantation used in complex function theory. In \cite{BaWa2_14} it was applied to some shape optimization problems involving Robin eigenvalues. For convenience we shortly review some of the principal properties. For his mehod we need the Green's function with Dirichlet boundary condition
\begin{eqnarray}\label{Green}
G_\Omega(x,y)= \gamma( S(|x-y|)-H(x,y)),
\end{eqnarray}
where
\begin{align}
\gamma = 
\begin{cases} 
\frac{1}{2\pi}&  \tx{if} n=2\\
\frac{1}{(n-2)|\p B_1|}& \tx{if} n>2
\end{cases}
&\qquad\tx{and} \qquad S(t)= \begin{cases}
-\ln(t) &\tx{if} n=2\\
t^{2-n} & \tx{if} n>2.
\end{cases}
\end{align}
For fixed $y\in \Omega$ the funcion $H(\cdot,y)$ is harmonic.
\begin{definition} The harmonic radius  at a point $y\in \Omega$ is given by
$$
r(y)= \begin{cases}
e^{-H(y,y)} &\tx{if} n=2,\\
H(y,y)^{-\frac{1}{n-2}}&\tx{if} n>2.
\end{cases}
$$
\end{definition}
The harmonic radius vanishes on the boundary $\p\Omega$ and takes its maximum $r_\Omega$ at the harmonic center $y_h$.  It satisfies the isoperimetric inequality (\cite{He69},\cite{BaFl96})
\begin{eqnarray}\label{harmonic radius}
|B_{r_\Omega}|\leq |\Omega|.
\end{eqnarray}
To illustrate the size of the harmonic radius we note that $r_{\Omega}$ is estimated from below by the inner radius $r_i(\Omega)$ and from above by the outer radius $r_{o}(\Omega)$ of the domain:
\begin{eqnarray*}
r_{i}(\Omega)\leq r_{\Omega}\leq r_{o}(\Omega).
\end{eqnarray*}
Note that $G_{B_R}(x,0)$ is a monotone function in $r=\vert x\vert$. Consider any radial function $\phi:B_{r_\Omega}\to \mathbb{R}$ thus $\phi(x)=\phi(r)$. Then there exists a function $\omega:\rz\to\rz$ such that
\begin{eqnarray*}
\phi(x)=\omega(G_{B_{r_\Omega}}(x,0)).
\end{eqnarray*}
To $\phi(x)$ we associate the transplanted function $U:\Omega \to \mathbb{R}$  defined  by $U(x)=\omega(G_\Omega(x,y_h))$. Then for any positive function $f(s)$, the following inequalities hold true
\begin{align}\label{harmonic1}
\int_\Omega |\nabla U|^2\:dx &= \int_{B_{r_\Omega}} \vert\nabla \phi|^2\:dx\\
\int_\Omega f(U)\:dx &\geq \int_{B_{r_\Omega}} f(\phi)\:dx.\label{harmonic2}\\
\int_\Omega f(U)\:dx &\leq \gamma^n\int_{B_{r_\Omega}} f(\phi)\:dx, \label{harmonic3}
\end{align}
where
\begin{eqnarray*}
\gamma = \left(\frac{|\Omega|}{|B_{r_\Omega}|}\right)^{\frac{1}{n}}.
\end{eqnarray*}
 For a proof see \cite{He69} or \cite{BaFl96} and in particular \cite{BaWa2_14} for a proof of \eqref{harmonic3}.
 The following observation will be useful in the sequel.
 \begin{remark}\label{iso}
 Since $U$ is constant on $\partial\Omega$ ($U=U(\partial\Omega)$) and since 
 $\phi$ is radial we deduce
 \begin{eqnarray*}
 \int\limits_{\partial\Omega}U^2\:dS
 =
 U^2(\partial\Omega)\:\vert\partial\Omega\vert=\phi^2(r_{\Omega})\:\vert\partial B_{r_{\Omega}}\vert\:\frac{\vert\partial\Omega\vert}{\vert\partial B_{r_{\Omega}}\vert}
 =
 \frac{\vert\partial\Omega\vert}{\vert\partial B_{r_{\Omega}}\vert}\:\int\limits_{\partial B_{r_{\Omega}}}\phi^2\:dS.
\end{eqnarray*}
Since by \eqref{harmonic radius}, $|B_{r_\Omega}|\leq |\Omega|$ the isoperimetric inequality implies that $\frac{|\p \Omega|}{|\p B_{r_\Omega}|}\geq1$.
\end{remark}
Consider first the case $\lambda^\sigma_1(\Omega)$ with $\sigma<\sigma_0=\frac{|\Omega|}{|\p \Omega|}.$
Let $u$ be a positive normalized radial eigenfunction of Problem \eqref{resolvent} in $B_{r_\Omega}$ with $\sigma$ replaced by $\sigma'$, corresponding to the eigenvalue $\lambda_{1}^{\sigma'}(B_{r_\Omega})$. Here  
\begin{eqnarray*}
\sigma'=\sigma\: \frac{\vert\partial\Omega\vert}{\vert\partial B_{r_{\Omega}}\vert}\leq \sigma.
\end{eqnarray*}
Since 
\begin{align}\label{sigma}
\sigma'= -|\sigma|\: \frac{\vert\partial\Omega\vert}{\vert\partial B_{r_{\Omega}}\vert}<\sigma_0(\Omega) \: \frac{\vert\partial\Omega\vert}{\vert\partial B_{r_{\Omega}}\vert}=-\frac{|\Omega|}{|\p B_{r_\Omega}|}\leq \sigma_0(B_{r_\Omega}),
\end{align}
$u$ is of constant sign.
Then the transplanted function $U$ of $u$ in $\Omega$ satisfies $\int_\Omega |\nabla U|^2\:dx=1$. By \eqref{bawa1} we have
\begin{align}
\frac{1}{\lambda_1(\Omega)}\geq \int_\Omega U^2\:dx +\sigma \oint_{\p\Omega} U^2\:dS.
\end{align}
Taking into account \eqref{harmonic2} and Remark \ref{iso} we get
\begin{eqnarray*}
\frac{1}{\lambda^{\sigma}_1(\Omega)}
\geq
\int\limits_{B_{r_{\Omega}}}u^2\:dx-\vert\sigma\vert\: \frac{\vert\partial\Omega\vert}{\vert\partial B_{r_{\Omega}}\vert}\:\int\limits_{\partial B_{r_{\Omega}}}u^2\:dS.
\end{eqnarray*}
By \eqref{sigma} the right-hand side is positive and is equal to $1/\lambda_1^{\sigma'}(B_{r_\Omega})$. Consequently
\begin{eqnarray*}
0\leq\lambda^{\sigma}_1(\Omega)
\leq
\lambda^{\sigma'}_1(B_{r_{\Omega}}).
\end{eqnarray*}
\medskip

Consider now the case $\sigma_0<\sigma<0$. Define
$$
\sigma''= \sigma \frac{|\p\Omega||B_{r_\Omega}|}{|\Omega||\p B_{r_\Omega}|}>-\frac{|B_{r_\Omega}|}{|\p B_{r_\Omega}}=\sigma_0(B_{r_\Omega}).
$$
Let $u$ be a positive normalized radial eigenfunction of problem \eqref{resolvent} in $B_{r_\Omega}$  corresponding to the eigenvalue $\lambda_{1}^{\sigma''}(B_{r_\Omega})$. Let $U$ be the transplanted function of $u$ in $\Omega$. Then by \eqref{bawa2} we get
\begin{eqnarray*}
\frac{1}{\lambda^\sigma_{-1}(\Omega)}\leq \int\limits_{\Omega}U^2\:dx-\vert\sigma\vert\int\limits_{\partial\Omega}U^2\:dS.
\end{eqnarray*}
We apply \eqref{harmonic3} to the first integral in the denominator and again Remark \ref{iso} to the second.
\begin{eqnarray*}
\frac{1}{\lambda^\sigma_{-1}(\Omega)}&\leq&\gamma^n\int\limits_{B_{r_{\Omega}}}u^2\:dx-\vert\sigma\vert\:\frac{\vert\partial\Omega\vert}{\vert\partial B_{r_{\Omega}}\vert}\:\int\limits_{\partial\Omega}u^2\:dS\\
&=&
\gamma^n\left(\int\limits_{B_{r_{\Omega}}}u^2\:dx+\sigma''\:\int\limits_{\partial\Omega}u^2\:dS\right).
\end{eqnarray*}
Thus
\begin{eqnarray*}
0>\lambda^\sigma_{-1}(\Omega)
\geq
\frac{1}{\gamma^n}\lambda^{\sigma''}_{-1}(B_{r_{\Omega}})
\end{eqnarray*}
We may rewrite this inequality as
\begin{eqnarray*}
\vert\Omega\vert\:\lambda^\sigma_{-1}(\Omega)
\geq
\vert B_{r_{\Omega}}\vert\:\lambda^{\sigma''}_{-1}(B_{r_{\Omega}}).
\end{eqnarray*}
This proves the following theorem.
\begin{theorem}\label{final}
Let $\Omega$ be any domain for which the trace operator $H^{1,2}(\Omega)\to L^2(\partial\Omega)$ is well defined. Let $\lambda^\sigma_{\pm1}(\Omega)$ be the first positive (negative) eigenvalue of \eqref{resolvent} given by \eqref{eigenchar}. Let $r_{\Omega}$ be the harmonic radius of 
$\Omega$. Then the following optimality result holds.
\begin{itemize}
\item[1)] In the case $\sigma<\sigma_0(\Omega)<0$ we have $0\leq\lambda^\sigma_{1}(\Omega)
\leq
\lambda^{\sigma'}_{1}(B_{r_{\Omega}})$,

\item[2)] In the case $\sigma_0(\Omega)<\sigma<0$ we have 
$0
>\vert\Omega\vert\:\lambda^\sigma_{-1}(\Omega)
\geq
\vert B_{r_{\Omega}}\vert\:\lambda^{\sigma''}_{-1}(B_{r_{\Omega}})$.
\end{itemize}
Equality holds in both cases if and only if $\Omega$ is a ball.
\end{theorem}
Since $\sigma'\leq\sigma$, Lemma \ref{sigmamon} gives
\begin{eqnarray*}
\lambda_{1}^{\sigma'}(B_{r_{\Omega}})\geq\lambda_{1}^{\sigma}(B_{r_{\Omega}}).
\end{eqnarray*}
In general $\sigma''$ and $\sigma$ are not comparable.
\begin{remark}\label{bawa3}
It is interesting to compare 1) in Theorem \ref{final}  with \eqref{glopt1} in Lemma \ref{domaindep} (i). We get the following  two sided bounds.
\newline
If $B_R$ is a ball of equal volume with $\Omega$ and if $\sigma<\sigma_0(B_R)<0$ then
\begin{eqnarray*}
\lambda^\sigma_{1}(B_R)\leq \lambda^\sigma_{1}(\Omega)
\leq
\lambda^{\sigma'}_{1}(B_{r_{\Omega}}).
\end{eqnarray*}
Equality holds for the ball.
\end{remark}
\section{Domain dependence}
\subsection{Small perturbations of a given domain}
We are interested in deriving optimality conditions for for the domain functionals 
$\lambda_{\pm 1}(\Omega)$. Contrary to the results in \cite{Ba08} (see \eqref{glopt1} and \eqref{glopt2}) these results will be local. We first decribe the general setting.
\newline
\newline
Consider a family of domains $(\Omega_t)_t$. The parameter $t$ varies in some open interval $(-t_0,t_0)$ where $t_0>0$ is prescribed. With this notation we set
$\Omega_0:=\Omega$. The family is given by the following construction. Let
\begin{eqnarray*}
\Phi_{t}:\Omega\to\Omega_{t}:=\Phi_{t}(\Omega)\qquad \qquad y:=\Phi_{t}(x)=x+tv(x)+\frac{t^2}{2}w(x)+o(t^2)
\end{eqnarray*}
be a smooth family, where  $v$ and $w$ are vector fields such that
\begin{equation*}
v,w:\overline{\Omega}\to\R\quad\hbox{are in}\: C^1(\overline{\Omega}).
\end{equation*}
Note that for $t_0>0$ small enough $(\Phi_{t})_{\vert t\vert<t_0}$ is a family of diffeomorphism. This restricts $t_0$ and defines the notion of "small perturbation of $\Omega$".
\newline
\newline
The volume of $\Omega_t$ is given by
\begin{eqnarray*}
\vert\Omega_t\vert=\int_{\Omega} J(t)\:dx
\end{eqnarray*}
where $J(t)$ is the Jacobian determinant corresponding to the transformation $\Phi_t$. The Jacobian matrix corresponding to this transformation is up to second order terms of the form
\begin{eqnarray*}
I + t D_v +\frac{t^2}{2}D_w, \tx{where} (D_v)_{ij}=\partial_j v_i \tx{and} \partial_j=\partial/ \partial x_j.
\end{eqnarray*}
By Jacobi's formula we have for small $t$
\begin{eqnarray}\label{jform}
J(t)&:=&\tx{det}(I+tD_v+\frac{t^2}{2}D_w)\\
\nonumber&=&1+t\:\Div v+\frac{t^2}{2}\left((\Div v)^2-D_v : D_v+\Div\: w\right)+o(t^2).
\end{eqnarray}
Here we used the notation
\begin{eqnarray*}
D_v : D_v:=\partial_{i}v_j\partial_{j}v_{i}.
\end{eqnarray*}
Hence
\begin{eqnarray}\label{volexp}
\vert\Omega_t\vert&=& \int_{B_R} J(t)\:dx=\vert\Omega\vert +t\int_{\Omega} \tx{div} v\:dx\\
\nonumber&&+ \frac{t^2}{2}\int_{\Omega}((\Div v)^2-D_v : D_v+\Div\:w)\:dx+o(t^2).
\end{eqnarray}
For the first variation we have only to require that $y$ is volume preserving of the first order, that is
\begin{eqnarray}\label{volume1}
\frac{d}{dt}\vert\Omega_t\vert\Big\vert_{t=0}=\int\limits_{\Omega} \Div v\:dx =\int\limits_{\partial \Omega}(v\cdot \nu)\:dS=0.
\end{eqnarray}
We also consider perturbations which, in addition to the condition \eqref{volume1}, satisfy the volume preservation of the second order, namely
\begin{eqnarray}\label{volume2}
\frac{d^2}{dt^2}\vert\Omega_t\vert\Big\vert_{t=0}=\int_{\Omega}((\Div v)^2-D_v : D_v+\Div\:w)\:dx=0.
\end{eqnarray}
In addition we consider perturbations which preserve the surface area up to second order. 
Instead of \eqref{volexp} we then use the expansion 
\begin{eqnarray*}
\vert\partial\Omega_{t}\vert=\int\limits_{\partial\Omega}m(t)\:dS
=
\vert\partial\Omega\vert
+
t\int\limits_{\partial\Omega}\dot{m}(0)\:dS
+
\frac{t^2}{2}\int\limits_{\partial\Omega}\ddot{m}(0)\:dS+o(t^2)
\end{eqnarray*}
From this we can derive first and second order conditions. 
In \cite{BaWa14} the following first order condition
\begin{eqnarray}\label{area1}
\int\limits_{\partial\Omega}(v\cdot\nu) H_{\partial\Omega}\:dS=0
\end{eqnarray}
and second order condition
\begin{eqnarray}\label{area2}
\int\limits_{\partial\Omega}F(\nabla^{*}v,v)\:dS
+
(n-1)\int\limits_{\partial\Omega}(w\cdot\nu) H_{\partial\Omega}\:dS
=0
\end{eqnarray}
were derived. Here $F(\nabla^{*}v,v)$ is a known scalar function of the tangential derivative $\nabla^{*}v$ of $v$ (see e.g. formula (2.20) in \cite{BaWa14}). Moreover, $H_{\partial\Omega}$ denotes the mean curvature of $\partial\Omega$. In particular for the ball $B_R$ this condition reads as
\begin{eqnarray}\label{area2ball}
\qquad\ddot{S}(0):=\int\limits_{\partial B_R}\vert\nabla^{*}(v\cdot\nu)\vert^2\:dS
-
\frac{n-1}{R^2}\int\limits_{\partial B_R}(v\cdot\nu)^2\:dS
+
\frac{n-1}{R}\frac{d^2}{dt^2}\vert\Omega_t\vert\Big\vert_{t=0}
=0,
\end{eqnarray}
where the last term in the sum is computed in \eqref{volume2} (see also Lemma 2 in \cite{BaWa14}).
For later use we set
\begin{eqnarray*}
\ddot{S}_{0}(0):=\int\limits_{\partial B_R}\vert\nabla^{*}(v\cdot\nu)\vert^2\:dS
-
\frac{n-1}{R^2}\int\limits_{\partial B_R}(v\cdot\nu)^2\:dS
\end{eqnarray*}
and
\begin{eqnarray*}
\ddot{V}(0):=\frac{d^2}{dt^2}\vert\Omega_t\vert\Big\vert_{t=0}.
\end{eqnarray*}
Thus 
\begin{eqnarray}\label{areanot}
\ddot{S}(0)=\ddot{S}_{0}(0)+\frac{n-1}{R}\ddot{V}(0).
\end{eqnarray}
Note that for volume preserving perturbations, $\ddot{S}(0)(=\ddot{S}_{0}(0))$ describes the 
isoperimetric defect and is strictly positive (see also Section 7 in \cite{BaWa14}).
\subsection{First and second domain variation}
\subsection{The first variation and monotonicity}
Let $(\Omega_t)_t$ be a smooth family of small perturbations of $\Omega$ as described in the previous subSection. In particular they will be either volume preserving in the sense of \eqref{volume1} and \eqref{volume2} or area preserving in the sense of  \eqref{area1} and \eqref{area2}. For the moment we denote by $\lambda$ either of the two first eigenvalues $\lambda_{\pm 1}$. We denote by $u_t(x):=u(y(x);t)$ the solution of 
\begin{eqnarray}\label{BBRt}
\Delta u_t+\lambda(\Omega_t) u_t=0\quad\hbox{in}\:\:\Omega_t,\qquad\partial_{\nu_t}u_t=\lambda(\Omega_t)\:\sigma\:u_t\quad\hbox{in}\:\:\partial\Omega_t.
\end{eqnarray}
Here $\lambda(\Omega_t)$ has the representation
\begin{eqnarray}\label{lambdat}
\lambda(t):=\lambda(\Omega_t)=\frac{1}{\int\limits_{\Omega_t}u_t^2\:dy+\sigma\int\limits_{\partial \Omega_t}u_t^2\:dS_t},
\end{eqnarray}
where $u_t$ solves \eqref{BBRt}. Consequently the energy is
\begin{eqnarray}\label{cale1}
\E(t)=\int\limits_{\Omega_t}\vert\nabla u_t\vert^2\:dy-\lambda(t)\left(\int\limits_{\Omega_t}u_t^2\:dy+\sigma\int\limits_{\partial \Omega_t}u_t^2\:dS_t\right)\equiv 0 \quad\hbox{for all $t$}.
\end{eqnarray}
In \cite{BaWa14} and more detailed in \cite{BaWa14a} the first and second variation of $\E$ with respect to $t$ were computed. For the first variation we obtained (see (4.1)) in \cite{BaWa14})
\begin{eqnarray}\label{ender}
0=\dot{\E}(0)&=&\int\limits_{\partial\Omega}(v\cdot\nu)
\left\{\vert\nabla u\vert^2-\lambda u^2-2\lambda^2\sigma^2u^2-\lambda\sigma(n-1)H_{\partial\Omega}u^2\right\}\:dS\\
&&\nonumber-\dot{\lambda}(0)\left(\int\limits_{\Omega}u^2\:dx+\sigma\int\limits_{\partial \Omega}u^2\:dS\right).
\end{eqnarray}
\begin{remark}\label{smoothness}
The differentiability of $\lambda(t)$ in $t=0$ is not automatic. In fact, the eigenvalues 
$\lambda=\lambda_{\pm 1}$ are differentiable in $t=0$ if $\lambda_{\pm 1}$ is simple (see e.g.  \cite[IV, Sec. 3.5]{Ka}). 
As we know from Lemma \ref{sigma0} (i) and (ii) in Section 2, this is true for our choice of $\sigma$.
\end{remark}
The condition $\dot{\lambda}(0)=0$ and \eqref{ender} gives the necessary condition
\begin{eqnarray*}
\int\limits_{\partial\Omega}(v\cdot\nu)
\left\{\vert\nabla u\vert^2-\lambda u^2-2\lambda^2\sigma^2u^2-\lambda\sigma(n-1)H_{\partial\Omega}u^2\right\}\:dS=0.
\end{eqnarray*}
In the case of volume preserving perturbation, \eqref{volume1} implies
\begin{eqnarray}\label{necvol}
\vert\nabla u\vert^2-\lambda u^2-2\lambda^2\sigma^2u^2-\lambda\sigma(n-1)H_{\partial\Omega}u^2=const.\qquad\hbox{on}\quad\partial\Omega.
\end{eqnarray}
This is a special case of Theorem 1 in \cite{BaWa14}. In the case of surface area preserving perturbations we apply \eqref{area1} and obtain
\begin{eqnarray}\label{necarea}
\vert\nabla u\vert^2-\lambda u^2-2\lambda^2\sigma^2u^2-\lambda\sigma(n-1)H_{\partial\Omega}u^2=const. H_{\partial\Omega}\qquad\hbox{on}\quad\partial\Omega
\end{eqnarray}
as a necessary condition for any critical point of $\lambda(\Omega)$. It is an open question whether \eqref{necvol} or \eqref{necarea} implies, that $\Omega$ can only be a ball.
\newline
\newline
From now on let $\Omega=B_R$ . To ensure differentiability of $\lambda$ in $t=0$ we consider 
the cases for
\begin{eqnarray}\label{keyass1}
\lambda =\lambda_{1}\qquad\hbox{if}\qquad\sigma<\sigma_{0}(B_R)
\end{eqnarray}
and
\begin{eqnarray}\label{keyass2}
\lambda =\lambda_{-1}\qquad\hbox{if}\qquad\sigma_{0}(B_R)<\sigma<0.
\end{eqnarray}
Then \eqref{necvol} and \eqref{necarea} are satisfied since the corresponding eigenfunctions are radial. Hence $\dot{\lambda}(0)=0$, i.e. the ball is a critical domain. 
\newline
\newline
Formula \eqref{ender} implies monotonicity of $\lambda_{\pm 1}$ for nearly spherical domains with respect to volume increasing (decreasing) perturbations. Indeed we rewrite \eqref{ender} as
\begin{eqnarray*}
a(u,u)\dot{\lambda}(0)=\int\limits_{ \partial B_R}\left(u_{r}^2(R)-\lambda u^2(R)-2\lambda^2\sigma^2u^2(R)-\lambda\sigma\frac{n-1}{R}u^2(R)\right)(v\cdot\nu)\:dS.
\end{eqnarray*}
Then we use the boundary condition $u_{r}^2(R)=\lambda^2\sigma^2u^2(R)$ and obtain
\begin{eqnarray}\label{lambdamon}
a(u,u)\dot{\lambda}(0)=-u(R)k(R)\int\limits_{ \partial B_R}(v\cdot\nu)\:dS,
\end{eqnarray}
where 
\begin{eqnarray}\label{eigenkr}
k(R):= \lambda\:u(R)\left(1+\frac{(n-1)\:\sigma}{R}+\lambda\:\sigma^2\right).
\end{eqnarray}
Next we determine the sign of $k(R)$. In this we modify the proof of Lemma 3 in \cite{BaWa2_14}. For the sake of completeness we give the details.
\begin{lemma}\label{kr}
Let $k(R)$ be given by \eqref{eigenkr} and let $u(r)$ be the positive radial function in the case $\lambda=\lambda_{1}$ 
or $\lambda=\lambda_{-1}$. Then we have
\begin{eqnarray*}
&&k(R)>0\qquad\qquad\hbox{if}\quad\lambda=\lambda_{1}
\\&&k(R)<0\qquad\qquad\hbox{if}\quad\lambda=\lambda_{-1}.
\end{eqnarray*}
\end{lemma}
{\bf{Proof}} In the radial case either eigenfunction satisfies the differentia equation
\begin{eqnarray*}
u_{rr}+\frac{n-1}{r}\:u_{r}+\lambda\:u(r)=0\quad\hbox{in}\:\:(0,R),\qquad u'(R)=\lambda\:\sigma\:u(R).
\end{eqnarray*}
We set $z=\frac{u_r}{u}$ and observe that
$$
\frac{dz}{dr} +z^2 +\frac{n-1}{r}z + \lambda=0  \tx{in} (0,R).
$$
At the endpoint
$$
\frac{dz}{dr}(R) + \lambda^2\:\sigma^2 +\frac{(n-1)}{R} \:\lambda\:\sigma + \lambda =0.
$$ 
We know that $z(0)=0$ and $z(R)= \lambda\:\sigma$. Note that 
\begin{eqnarray}\label{zinc}
z_r(0)= -\lambda.
\end{eqnarray}
We distinguish two cases.
\newline
\newline
The case $\lambda=\lambda_{-1}(B_R)$.
\newline
In that case we have (see also \eqref{zinc})
\begin{eqnarray}\label{zink2}
z(0)=0\:,\qquad z(R)=\lambda_{-1}\:\sigma>0\:,\qquad z_r(0)= -\lambda_{-1}>0.
\end{eqnarray}
Thus $z(r)$ increases near $0$. We again determine the sign of $z_r(R)$. If $z_r(R)\leq 0$ then because of \eqref{zink2} there exists a number $\rho \in (0,R)$ such that $z_r(\rho) =0$, $z(\rho)>0$ and $z_{rr}(\rho) \leq 0$. 
From the equation we get $z_{rr}(\rho) =\frac{n-1}{\rho^2} z(\rho)>0$ which is contradictory. Consequently
\begin{eqnarray*}
z_r(R)= -\left( \lambda_{-1}^2\:\sigma^2 +\frac{(n-1)}{R} \:\lambda_{-1}\:\sigma +\lambda_{-1} \right)>0.
\end{eqnarray*}
This also implies $k(R)<0$ in the case $\lambda=\lambda_{-1}(B_R)$.
\newline
\newline
The case $\lambda=\lambda_{1}(B_R)$.
\newline
We have (also from \eqref{zinc})
\begin{eqnarray}\label{zink1}
z(0)=0\:,\qquad z(R)=\lambda_{1}\:\sigma<0\:,\qquad z_r(0)= -\lambda_{1}<0.
\end{eqnarray}
Thus $z(r)$ decreases near $0$. We determine the sign of $z_r(R)$. If $z_r(R)\geq 0$ then because of \eqref{zink1} there exists a number $\rho \in (0,R)$ such that $z_r(\rho) =0$, $z(\rho)<0$ and $z_{rr}(\rho) \geq 0$. 
From the equation we get $z_{rr}(\rho) =\frac{n-1}{\rho^2} z(\rho)<0$ which leads to a contradiction. Consequently
\begin{eqnarray*}
z_r(R)= -\left( \lambda_{1}^2\:\sigma^2 +\frac{(n-1)}{R} \:\lambda_{1}\:\sigma +\lambda_{1} \right)<0.
\end{eqnarray*}
This implies $k(R)>0$ in the case $\lambda=\lambda_{1}(B_R)$.
\hfill $\square$
\newline
\newline
We easily prove the following lemma.
\begin{lemma}\label{lambdamon2}
The first derivative $\dot{\lambda}_{\pm 1}(0)$ satisfies the following sign condition
\begin{eqnarray*}
\dot{\lambda}_{\pm 1}(0) &<& 0\qquad\hbox{if} \quad \int\limits_{ \partial B_R}(v\cdot\nu)\:dS>0,\\
\dot{\lambda}_{\pm 1}(0) &>& 0\qquad\hbox{if} \quad \int\limits_{ \partial B_R}(v\cdot\nu)\:dS<0.
\end{eqnarray*}
\end{lemma}
{\bf{Proof}} This follows directly from \eqref{lambdamon} and Lemma \ref{kr}. Indeed it is sufficient to recall that $a(u,u)>0$ for $\lambda=\lambda_{1}(B_R)$ and $a(u,u)<0$ for $\lambda=\lambda_{-1}(B_R)$. Also note that in either case $u$ is positive.
\hfill $\square$
\subsubsection{The second variation}
We are interested in extremality properties of the ball. We set $u'=\partial_{t}u_t(x)\vert_{t=0}$. This quantity is called "shape derivative" and plays a crucial role in determining the sign of the second domain variation of $\lambda$. In a first step we derive an equation for $u'$. This follows the technique in \cite{BaWa2_14}. Another good reference, with a slightly different approach, is also the book of A. Henrot an M. Pierre \cite{HePi05}. As a result we consider the following boundary value problem.
\begin{eqnarray}\label{eigenprim}
\Delta u'+\lambda\:u'=0\quad\hbox{in}\:\: B_R\qquad\partial_{\nu}u'-\lambda\:\sigma\:u'=k(R)(v\cdot\nu)\quad\hbox{in}\:\:\partial B_R
\end{eqnarray}
where $k(R)$ is given in \eqref{eigenkr}.
\medskip

To \eqref{eigenprim} we associate the quadratic form
\begin{eqnarray*}
Q(u'):=\int\limits_{B_R}\vert\nabla u'\vert^2\:dx-\lambda\:\int\limits_{B_R}u'^2\:dx-\lambda\:\sigma\: \int\limits_{\partial B_R}u'^2\:dS.
\end{eqnarray*}
We now turn to the computation of $\ddot{\lambda}(0)$. If we write \eqref{cale1} as
\begin{eqnarray*}
\E(t)=\F_{1}(t)-\lambda(t)\F_{2}(t).
\end{eqnarray*}
If we use $\dot{\lambda}(0)=0$, we obtain the formula
\begin{eqnarray*}
\ddot{\lambda}(0)=\frac{\ddot{\F}_{1}(0)-\lambda(0)\ddot{\F}_{2}(0)}{\F_{2}(0)}.
\end{eqnarray*}
In the case $\Omega=B_R$, we repeat computations as done in \cite{BaWa14} (Section 7). However 
we don't assume any more, that the perturbations are volume preserving. This leads to the following
modified formula (note that $\F_2(0)=a(u,u)$)
\begin{eqnarray}\label{lambda2}
\ddot{\lambda}(0)a(u,u)
&=&
-2Q(u')
-
\lambda\sigma \:u^2(R)\:\ddot{S}_{0}(0)
-
k(R)\:u(R)\ddot{V}(0)\\
\nonumber&&-
2\lambda\sigma k(R)u(R)\int\limits_{\partial B_R}(v\cdot\nu)^2\:dS.
\end{eqnarray}
In view of \eqref{eigenchar}, \eqref{keyass1} and \eqref{keyass2} the ball $B_R$ is a local minimizer for $\lambda$ if
\begin{itemize}
\item[(1)] $\lambda=\lambda_{1}>0$, $\sigma<\sigma_{0}(R)<0$ and
\begin{eqnarray}\label{necess1}
&&-2Q(u')
-
\lambda_{1}\sigma \:u^2(R)\:\ddot{S}_{0}(0)
-
k(R)\:u(R)\ddot{V}(0)\\
&&\nonumber\qquad-
2\lambda_{1}\sigma k(R)u(R)\int\limits_{\partial B_R}(v\cdot\nu)^2\:dS>0;
\end{eqnarray}
\item[(2)] $\lambda=\lambda_{-1}<0$, $\sigma_{0}(R)<\sigma<0$ and
\begin{eqnarray}\label{necess2}
&&-2Q(u')
-
\lambda_{-1}\sigma \:u^2(R)\:\ddot{S}_{0}(0)
-
k(R)\:u(R)\ddot{V}(0)\\
&&\nonumber\qquad
-
2\lambda_{-1}\sigma k(R)u(R)\int\limits_{\partial B_R}(v\cdot\nu)^2\:dS<0.
\end{eqnarray}
\end{itemize}
In the next section we will discuss the sign of $\ddot{\lambda}(0)$.  
\subsubsection{The sign of the second variation}
We consider the following Steklov eigenvalue problem
\begin{align}\label{steklovg}
\Delta \phi +\lambda\:\phi=0 \tx{in} B_R,\\
\nonumber \partial_{\nu} \phi -\lambda\:\sigma\:\phi =\mu\: \phi \tx{on} \p B_R.
\end{align}
There exists an infinite number of eigenvalues
\begin{eqnarray*}
\mu_1<\mu_2\leq \mu_3\leq ... \lim_{i\to \infty} \mu_i=\infty.
\end{eqnarray*}
\begin{remark}\label{spec}
Note that for $\lambda=\lambda_{\pm 1}$ we have an eigenvalue $\mu=0$. Indeed, the case 
$\mu=0$ corresponds to the case where $\phi=u_{\pm 1}$. For $\sigma<\sigma_{0}(B_R)$ (resp.  $\sigma_{0}(B_R)<\sigma<0$) the eigenvalue $\lambda_{1}$ (resp. $\lambda_{-1}$) is simple and the eigenfunction $u_{1}$ (resp. $u_{-1}$) is of constant sign. Thus $0=\mu=\mu_1$. As a consequence the spectrum consists only of non-negative eigenvalues.
\end{remark}
There exists a complete system of eigenfunctions $\{\phi_i\}_{i\geq 1}$ such that
\begin{eqnarray*}
\int\limits_{\partial\Omega}\phi_{i}\phi_{j}\:dS=\delta_{ij}.
\end{eqnarray*}
Similarly to the dicussion in Section 7 in \cite{BaWa14} we get the representation 
\begin{eqnarray*}
u'=\sum_{i=2}^{\infty}c_{i}\:\phi_{i}\qquad(v\cdot\nu)=\sum_{i=2}^{\infty}b_{i}\phi_i
\end{eqnarray*}
for the solution $u'$ of \eqref{eigenprim} and the perturbation $v\cdot\nu$. 
Note that by \eqref{volume1} or \eqref{area1} we have $c_1=b_1=0$.
\newline
\newline
It is easy to check that
\begin{eqnarray}\label{no2}
Q(u')=\sum_{i=2}^{\infty}c_{i}^2\mu_{i}.
\end{eqnarray}
The boundary condition in \eqref{eigenprim} implies
\begin{eqnarray}\label{no3}
b_{i}=\frac{c_{i}\:\mu_{i}}{k(R)}\qquad\hbox{thus}\qquad \int\limits_{\partial B_R}(v\cdot\nu)^2\:dS_R=\sum_{i=2}^{\infty}\frac{c_i^2\mu_i^2}{k^2(R)},
\end{eqnarray}
where $k(R)$ is defined in \eqref{eigenkr}. We will also use the estimate
\begin{eqnarray}\label{no1}
\ddot{S}_{0}(0)\geq \frac{n+1}{R^2}\int\limits_{\partial B_R}(v\cdot\nu)^2\:dS=\frac{n+1}{k^2(R)R^2}\sum_{i=2}^{\infty}c_i^2\mu_i^2.
\end{eqnarray}
This was shown in the derivation of (7.13) in \cite{BaWa14} and holds equally in the case \eqref{volume1} or \eqref{area1}.
\medskip

\noindent{\bf{Volume preserving perturbations}}
\medskip

\noindent
In this case we have $\ddot{V}(0)=0$.  We first consider $\lambda_{1}$. According to \eqref{necess1}, a necessary condition for the ball to be a minimizer is
\begin{eqnarray}\label{secpoei}
-2Q(u')-\lambda_{1}\sigma u^2(R)\ddot{S}_{0}(0)-2\lambda\sigma k(R) u(R)\int\limits_{\partial B_R}(v\cdot\nu)^2\:dS
>0.
\end{eqnarray}
Since $\lambda_{1}>0$ and $\sigma<0$ we can apply \eqref{no1}, \eqref{no2} and \eqref{no3}. Thus it is sufficient to show that
\begin{eqnarray*}
2\sum_{i=1}^{\infty}c_i^2\:\mu_i^2\left\{-\frac{1}{\mu_{2}} + \frac{(n+1)\lambda_{1}\vert\sigma\vert\:u^2(R)}{2R^2k^2(R)}+\frac{\lambda_{1}\vert\sigma\vert u(R)}{k(R)}\right\}>0,
\end{eqnarray*}
where $\mu_2$ is the first positive eigenvalue of the Steklov eigenvalue problem \eqref{steklovg}. Such an expression also appeared in \cite{BaWa14} Section 7.2.2 subsection 2 where the second domain variation of a Robin eigenvalue problem was considered. It was shown that 
\begin{eqnarray*}
L:=\mu_2-\alpha+\frac{n-1}{R}-\frac{\lambda}{\alpha}=0
\end{eqnarray*}
If we set $\lambda=\lambda_{1}$, $\alpha=-\lambda_{1}\sigma$ and take into account \eqref{eigenkr} this is equivalent to
\begin{eqnarray*}
\frac{\lambda_{1}\vert\sigma\vert\:u(R)}{k(R)} -\frac{1}{\mu_{2}}=0.
\end{eqnarray*}
This proves \eqref{secpoei}.
\begin{theorem}\label{optimality1}
For some given ball $B_R$ let $\sigma_0(B_R)=-\frac{R}{n}$. Let $\lambda_{1}$ be given by 
\eqref{eigenchar}. Then the following optimality result holds: If $\sigma<\sigma_{0}(B_R)$, then
among all smooth domains of equal volume the ball $B_R$ is a local minimizer for 
$\lambda_{1}$. If we exclude translations and rotations of $B_R$ then it is a strict local minimizer.
\end{theorem}
The case $\lambda=\lambda_{-1}(B_R)<0$ and $\sigma_{0}(R)<\sigma<0$ is similar. 
In that case $k(R)<0$ by Lemma \ref{kr}. Since $\ddot{V}(0)=0$, we deduce from \eqref{necess2} the following necessary condition for the ball to be a minimizer:
\begin{eqnarray}\label{secneei}
2Q(u')
+
\lambda_{-1}\sigma \:u^2(R)\:\ddot{S}_{0}(0)
+
2\lambda_{-1}\sigma k(R)u(R)\int\limits_{\partial B_R}(v\cdot\nu)^2\:dS>0.
\end{eqnarray}
Note that in this case $\lambda_{-1}\sigma>0$, thus we can apply \eqref{no1} again. As a consequence \eqref{secneei} holds, if 
\begin{eqnarray*}
2Q(u')
+
\left(\lambda_{-1}\sigma \:u^2(R)\frac{n+1}{R^2}
+
2\lambda_{-1}\sigma k(R)u(R)\right)\int\limits_{\partial B_R}(v\cdot\nu)^2\:dS>0.
\end{eqnarray*}
In view of \eqref{no2} and \eqref{no3} this is equivalent to
\begin{eqnarray*}
2\sum_{i=1}^{\infty}c_i^2\:\mu_i^2\left\{ \frac{1}{\mu_{2}}+\frac{(n+1)\lambda_{-1}\sigma\:u^2(R)}{2R^2 k^2(R)} +\frac{\lambda_{-1}\sigma u(R)}{k(R)} \right\}>0.
\end{eqnarray*}
Note that the first term in the sum is positive, the second is positive as well, while the third term is negative. Consequently it suffices to show
\begin{eqnarray*}
\frac{(n+1)\lambda_{-1}\sigma\:u^2(R)}{2R^2 k^2(R)} +\frac{\lambda_{-1}\sigma u(R)}{k(R)}>0.
\end{eqnarray*}
This is equivalent to
\begin{eqnarray}\label{nec1}
\frac{n+1}{R^2}+2\lambda_{-1}\left(1+\frac{n-1}{R}\sigma+\lambda_{-1}\sigma^2\right)>0.
\end{eqnarray}
Since 
\begin{eqnarray*}
-\frac{R}{n}=\sigma_{0}<\sigma<0
\end{eqnarray*}
we set $\sigma:=-\delta \frac{R}{n}$ for some $0<\delta<1$. Then \eqref{nec1} reads as
\begin{eqnarray}\label{nec2}
\frac{n+1}{R^2}+2\lambda_{-1}\left(1-\frac{n-1}{n}\delta\right)+2\lambda^2_{-1}\delta^2\frac{R^2}{n^2}>0.
\end{eqnarray}
Inequality \eqref{nec2} is a quadratic inequality in $\lambda_{-1}$. It is easy to check that both zeros of the quadratic exprssion are negative for all $0<\delta\leq 1$. This proves the minimality of the ball for $\lambda_{-1}$ for all $\sigma_{0}(B_R)<\sigma<0$.
\begin{theorem}\label{optimality2}
For some given ball $B_R$ let $\sigma_0(B_R)=-\frac{R}{n}$. Let $\lambda_{-1}$ be given by \eqref{eigenchar}. Then the following optimality result holds: If $\sigma_0(B_R)<\sigma<0$, then among all smooth domains of equal volume the ball $B_R$ is also a local minimizer for 
$\lambda_{-1}$. If we exclude translations and rotations of $B_R$ then it is a strict local minimizer.
\end{theorem}
\medskip

\noindent{\bf{Area preserving perturbations}}
\medskip

\noindent
In this case $\ddot{S}(0)$ and \eqref{areanot} gives
\begin{eqnarray}\label{volref}
\ddot{V}(0)=-\frac{R}{n-1}\:\ddot{S}_{0}(0)<0
\end{eqnarray}
We first consider $\lambda_{1}$. According to \eqref{necess1} and \eqref{volref}, a necessary condition for the ball to be a minimizer is
\begin{eqnarray}\label{asecpoei}
&&-2Q(u')-\lambda_{1}\sigma u^2(R)\ddot{S}_{0}(0)+  \frac{R}{n+1}k(R)u(R)\ddot{S}_{0}(0)
\\
&&\nonumber\qquad-2\lambda\sigma k(R) u(R)\int\limits_{\partial B_R}(v\cdot\nu)^2\:dS_R
>0.
\end{eqnarray}
If we compare this with \eqref{secpoei} we see that an extra positive term occurs. Since the arguments in \cite{BaWa14} Section 7.2.2 subSection 2 carry over to the area preserving case without any changes we proved the following theorem.
\begin{theorem}\label{optimality3}
For some given ball $B_R$ let $\sigma_0(B_R)=-\frac{R}{n}$. Let $\lambda_{1}$ be given by 
\eqref{eigenchar}. Then the following optimality result holds: If $\sigma<\sigma_{0}(B_R)$, then
among all smooth domains of equal area the ball $B_R$ is a local minimizer for 
$\lambda_{1}$. If we exclude translations and rotations of $B_R$ then it is a strict local minimizer.

\end{theorem}
The case $\lambda=\lambda_{-1}(B_R)<0$ and $\sigma_{0}(R)<\sigma<0$ is similar. 
In that case $k(R)<0$ by Lemma \ref{kr}. Since \eqref{volref} holds, we deduce from \eqref{necess2} the following necessary condition for the ball to be a minimizer:
\begin{eqnarray}\label{asecneei}
2Q(u')
+
\lambda_{-1}\sigma \:u^2(R)\:\ddot{S}_{0}(0)
+
2\lambda_{-1}\sigma k(R)u(R)\int\limits_{\partial B_R}(v\cdot\nu)^2\:dS>0.
\end{eqnarray}
Note that in this case $\lambda_{-1}\sigma>0$, thus we can apply \eqref{no1} again. As a consequence \eqref{asecneei} holds, if 
\begin{eqnarray*}
&& 2Q(u')
-
\frac{R}{n+1}k(R)u(R)\ddot{S}_{0}(0)\\
&&\nonumber\qquad
+
\left(\lambda_{-1}\sigma \:u^2(R)\frac{n+1}{R^2}
+
2\lambda_{-1}\sigma k(R)u(R)\right)\int\limits_{\partial B_R}(v\cdot\nu)^2\:dS>0.
\end{eqnarray*}
Again an additional positive term occurs. 
\begin{theorem}\label{optimality4}
For some given ball $B_R$ let $\sigma_0(B_R)=-\frac{R}{n}$. Let $\lambda_{-1}$ be given by \eqref{eigenchar}. Then the following optimality result holds: If $\sigma_0(B_R)<\sigma<0$, then among all smooth domains of equal area the ball $B_R$ is also a local minimizer for 
$\lambda_{-1}$. If we exclude translations and rotations of $B_R$ then it is a strict local minimizer.
\end{theorem}
\section{Open problems}
1. The variational characterization of $\lambda_{\pm 1}$ (see \eqref{eigenchar}) is also related to two inequalities, known as Friedrich's inequality and trace inequality. In fact, for $\sigma<\sigma_0<0$ we get
\begin{eqnarray}\label{fried}
\int\limits_{\Omega}u^2\:dx
\leq
\frac{1}{\lambda_{1}(\Omega)}\int\limits_{\Omega}\vert\nabla u\vert^2\:dx
+
\vert\sigma\vert\int\limits_{\partial\Omega}u^2\:dS \qquad\hbox{Friedrich's inequality}
\end{eqnarray}
and for $\sigma_0<\sigma<0$ we get
\begin{eqnarray}\label{trace}
\int\limits_{\partial\Omega}u^2\:dS
\leq
\frac{1}{\sigma\lambda_{-1}(\Omega)}\int\limits_{\Omega}\vert\nabla u\vert^2\:dx
+
\frac{1}{\vert\sigma\vert}\int\limits_{\Omega}u^2\:dx \qquad\hbox{trace inequality}.
\end{eqnarray}
The second inequality was also considered in \cite{Au14} where also the case of equality was analyzed. 
\newline
\newline
It is an interesting open problem to find explicit lower bounds for $\lambda_{1}(\Omega)$ and $\sigma\lambda_{-1}(\Omega)$. 
Note that the technique of harmonic transplantation gives upper bounds for these two quantities.
\newline
\newline
2. At least for $\lambda_{-1}(\Omega)$ it may be true that the ball of equal volume is only a local minimizer. There is no global result available at the moment. Therefore - motivated by Theorem \ref{final} 2) - it may also be interesting to ask if a quantity like $\vert\Omega\vert\lambda_{-1}(\Omega)$ has the ball of equal volume as a minimizer.  
\newline
\newline
{\bf Acknowledgement} This paper was written during a visit at the Newton Institute in Cambridge. Both authors would like to thank this Institute for the excellent working atmosphere. 

\end{document}